\documentclass[12pt]{article}
\usepackage{graphicx}
\usepackage{amsmath,amsthm,amssymb,enumerate}
\usepackage{euscript,mathrsfs}
\usepackage{dsfont}
\usepackage[left=2cm,right=2cm,top=3.5cm,bottom=3.5cm]{geometry}
\usepackage{color}
\catcode`\@=11 \@addtoreset{equation}{section}

\catcode`\@=12

\usepackage{accents}
\newcommand*{\tecka}[1]{%
  \accentset{\mbox{\large\bfseries .}}{#1}}

\allowdisplaybreaks

\newcommand\tr{\operatorname{tr}}

\newtheorem{Theorem}{Theorem}[section]
\newtheorem{Proposition}[Theorem]{Proposition}
\newtheorem{Lemma}[Theorem]{Lemma}
\newtheorem{Corollary}[Theorem]{Corollary}

\theoremstyle{definition}
\newtheorem{Definition}[Theorem]{Definition}

\newtheorem{Remark}[Theorem]{Remark}

\newcommand{\bTheorem}[1]{
\begin{Theorem} \label{T#1} }
\newcommand{\eT}{\end{Theorem}}

\newcommand{\bProposition}[1]{
\begin{Proposition} \label{P#1}}
\newcommand{\eP}{\end{Proposition}}
\newcommand{\be}{b_\ep}
\newcommand{\bLemma}[1]{
\begin{Lemma} \label{L#1} }
\newcommand{\eL}{\end{Lemma}}

\newcommand{\bCorollary}[1]{
\begin{Corollary} \label{C#1} }
\newcommand{\eC}{\end{Corollary}}

\newcommand{\bRemark}[1]{
\begin{Remark} \label{R#1} }
\newcommand{\eR}{\end{Remark}}

\newcommand{\bDefinition}[1]{
\begin{Definition} \label{D#1} }
\newcommand{\eD}{\end{Definition}}

\newcommand{\Del}{\Delta_x}

\newcommand{\bfj}{\mathbf{j}}

\newcommand{\bfphi}{\boldsymbol{\varphi}}

\newcommand{\bFormula}[1]{
\begin{equation} \label{#1}}
\newcommand{\eF}{\end{equation}}

\newcommand{\Ov}[1]{\overline{#1}}

\newcommand{\DC}{C^\infty_c}

\newcommand{\vr}{\varrho}
\newcommand{\vre}{\vr_\ep}

\newcommand{\vue}{{\vc v}_\ep}

\newcommand{\vu}{\vc{v}}
\newcommand{\vm}{\vc{m}}

\newcommand{\vc}[1]{{\bf #1}}

\newcommand{\Div}{{\rm div}_x}
\newcommand{\Grad}{\nabla_x}

\newcommand{\dx}{\,{\rm d} {x}}

\newcommand{\dt}{\,{\rm d} t }

\newcommand{\intO}[1]{\int_{\Omega} #1 \ \dx}

\newcommand{\vv}{\vc{v}}

\newcommand{\ep}{\varepsilon}

\definecolor{Cgrey}{rgb}{0.85,0.85,0.85}
\definecolor{Cblue}{rgb}{0.50,0.85,0.85}
\definecolor{Cred}{rgb}{1,0,0}
\definecolor{fancy}{rgb}{0.10,0.85,0.10}

\newcommand\Cbox[2]{%
    \newbox\contentbox%
    \newbox\bkgdbox%
    \setbox\contentbox\hbox to \hsize{%
        \vtop{
            \kern\columnsep
            \hbox to \hsize{%
                \kern\columnsep%
                \advance\hsize by -2\columnsep%
                \setlength{\textwidth}{\hsize}%
                \vbox{
                    \parskip=\baselineskip
                    \parindent=0bp
                    #2
                }%
                \kern\columnsep%
            }%
            \kern\columnsep%
        }%
    }%
    \setbox\bkgdbox\vbox{
        \color{#1}
        \hrule width  \wd\contentbox %
               height \ht\contentbox %
               depth  \dp\contentbox
        \color{black}
    }%
    \wd\bkgdbox=0bp%
    \vbox{\hbox to \hsize{\box\bkgdbox\box\contentbox}}%
    \vskip\baselineskip%
}

\begin{document}


\title{On a class of compressible viscoelastic rate-type fluids with stress--diffusion}

\author{Miroslav Bul\'{\i}\v{c}ek$^{**}$, Eduard Feireisl$^{*}$, Josef M\'{a}lek$^{**}$}




\maketitle

\bigskip

\centerline{$^{*}$ Institute of Mathematics, Czech Academy of Sciences}

\centerline{\v Zitn\' a 25, CZ-115 67 Prague 1, Czech Republic}

\centerline{}

\centerline{$^{**}$ Charles University, Faculty of Mathematics and Physics, Mathematical Institute}

\centerline{Sokolovsk\'a 83, CZ-186 75 Prague 8, Czech Republic}

\bigskip

\begin{abstract}

We develop a mathematical theory for a class of compressible viscoelastic rate-type fluids with stress diffusion. Our approach is based on the concepts used
in the nowadays standard theory of compressible Newtonian fluids as renormalization, effective viscous flux identity, compensated compactness. The presence
of the extra stress, however, requires substantial modification of these techniques, in particular, a new version of the effective viscous flux identity is derived.
With help of these tools, we show the existence of global--in--time weak solutions for any finite energy initial data.

\end{abstract}



\section{Introduction}
\label{I}

The objective of this paper is to develop a rigorous mathematical theory for a class of compressible viscoelastic rate-type fluid models with stress-diffusion. The governing equations are obtained using a thermodynamical approach based on the concept of evolving natural configuration that splits the total deformation process into a part that is purely elastic and another part that takes into account the  irreversible changes in the natural configuration by dissipation. The class of fluids considered in this study is such that the elastic response is described by a spherical strain (a scalar multiplier of the identity tensor). The model is derived in Appendix (carrying on the developments presented in \cite{bmps2018}, \cite{mp2016}, \cite{mtr2015}, \cite{mpss2018} and \cite{RASR}).

Recently, a similar class of models has been analyzed within the context of incompressible fluids in \cite{bmps2018}. In
the present study we establish long-time and large-data existence of weak solution to the evolutionary spatially periodic problem associated with the set of governing equations for compressible fluids. As the fluid is supposed to be compressible, our approach is based on the ideas borrowed from the theory of compressible \emph{Newtonian} fluids introduced by Lions \cite{LI4} and later elaborated in  \cite{EF70}. We mainly focus on the part of the theory based on the so--called Lions' identity
involving the crucial quantity - the effective viscous flux. The presence of the extra stress tensor in the present model destroys the original structure of \cite{LI4}, and, consequently, this part of the theory must by substantially modified.

The structure of the present paper is adapted to the mathematical theory - a detailed derivation of the model is given in Appendix. After a brief introduction of the problem in the preamble of Section \ref{M}, we state our main result. In Section \ref{AP}, we derive the necessary {\it a priori} bounds based on the energy balance. Sections \ref{PP}, \ref{S} form the heart of the paper. Here we derive a new version of Lions' identity for the effective viscous flux and use it for obtaining the necessary pressure bounds as well as weak sequential stability (compactness) of the family of weak solutions. Finally, in Section \ref{E}, we propose a suitable approximation scheme yielding the desired long time and large data existence result.

\section{Mathematical model, main result}
\label{M}

The mathematical model considered in the present paper is derived in Appendix. Beside the standard field - the fluid mass density $\vr = \vr(t,x)$ and the bulk velocity
$\vu = \vu(t,x)$ - the family of unknowns is supplemented by a new variable $b = b(t,x)$ that stands for the spherical part of the elastic strain, see Appendix.

For the sake of simplicity, we impose the periodic boundary conditions in $R^N$, $N=1,2,3$. Equivalently, we consider the equations on the spatial domain $\Omega$ identified with the flat torus,
\[
\Omega = \left( [-1,1]|_{ \{ -1, 1 \} } \right)^N,\ N=1,2,3.
\]
In fact, we focus on the physically relevant and technically most difficult case $N=3$. Note that substantial improvement of the theory can be expected for
$N = 1$, where mathematics of fluid flows of viscous compressible fluids is quite well understood, see e.g. the monograph by Antontsev~et~al.~\cite{AKM}.

\subsection{Field equations, constitutive relations}

The basic system of partial differential equations to be solved for $t \in (0,T)$, $x \in \Omega \subset R^3$ and the unknowns $\vr$, $\vu$, and $b$ reads:
\begin{equation} \label{M1}
\partial_t \vr + \Div (\vr \vu) = 0,
\end{equation}
\begin{equation} \label{M2}
\begin{split}
\partial_t (\vr \vu) + \Div (\vr \vu \otimes \vu) &+ \Grad \left( p_{{\rm fl}} (\vr) + p_{{\rm el}}(b)  + \frac{2}{3}
\sigma b \Del b \right) \\
&= \Div \mathbb{S} (\Grad \vu) - \sigma \Div \left( \Grad b \otimes \Grad b - \frac{1}{2} |\Grad b |^2 \mathbb{I} \right),
\end{split}
\end{equation}
\begin{equation} \label{M3}
\partial_t b + \vu \cdot \Grad b + \frac{1}{\nu}\left(e'(b)  - \sigma \Del b\right) = \frac{2}{3} b \Div \vu.
\end{equation}
Regarding the tensor $\mathbb{S}$ we assume that
\begin{equation} \label{M4}
\mathbb{S}(\Grad \vu) = \mu \left( \Grad \vu + (\Grad \vu)^{T} - \frac{2}{3} \Div \vu \mathbb{I} \right) + \lambda \Div \vu \mathbb{I},
\end{equation}
where $(\Grad \vu)^{T}$ denotes the transpose tensor to $\Grad \vu$. Moreover, we focus on the case where the quantities $\mu > 0$, $\lambda \geq 0$, $\nu > 0$, and $\sigma > 0$ are constants.

In addition, we suppose that $e : [0, \infty) \to [0, \infty]$ is a convex function, whereas
\begin{equation} \label{M5}
p_{\rm el} (b) = - e(b) - \frac{2}{3} b e'(b).
\end{equation}
The specific example we have in mind reads (see \eqref{aA12} below)
\begin{equation} \label{M6}
e(b) = a_1 (b - 1) - a_2 \log b ,
\end{equation}
but more general cases can be considered.

The problem is (formally) closed by imposing the initial conditions:
\begin{equation} \label{IC}
\vr(0, \cdot) = \vr_0, \ (\vr \vu)(0, \cdot) = \vm_0, \ b(0, \cdot) = b_0.
\end{equation}

\subsection{Weak solutions, main result}

Our goal is to establish the existence of global--in--time weak solutions to problem (\ref{M1}--\ref{IC}). In view of future applications, in particular to
problems concerning stability (cf. \cite{FeJiNo}), it is convenient to include the associated energy balance as a part of the definition of admissible weak solutions.
As shown in Section \ref{AP}, the total energy associated with the problem (\ref{M1}--\ref{M3})  takes the form  (see also \eqref{aA26} and \eqref{aA20c})
\begin{equation} \label{M7}
E(\vr, b, \Grad b,  \vu) = \frac{1}{2} \vr |\vu|^2 + \frac{\sigma}{2} |\Grad b |^2 + \Psi(\vr) +  e(b) ,\
\qquad \Psi ' (\vr) \vr - \Psi(\vr) = p_{\rm fl}(\vr),
\end{equation}
and, at least for any smooth spatially periodic solution, the energy balance
\[
\begin{split}
\frac{{\rm d}}{{\rm d}t} &\intO{ \left[ \frac{1}{2} \vr |\vu|^2 + \Psi(\vr) + \frac{\sigma}{2} |\Grad b|^2 +  e(b) \right] } \\
&+ \intO{ \mathbb{S}(\Grad \vu) : \Grad \vu } + \intO{ \frac{1}{\nu}\left( e'(b) - \sigma \Del b  \right)^2 } = 0.
\end{split}
\]
holds.

This motivates the following definition.

\begin{Definition} \label{D1}

The trio $[\vr, \vu, b]$ is called a \emph{dissipative weak solution} to problem (\ref{M1}--\ref{IC}) if:

\begin{itemize}
\item $\vr \geq 0$, $b > 0$ a.a. in $(0,T) \times \Omega$;
\item the integral identity
\[
\int_0^T \intO{ \left[ \vr \partial_t \varphi + \vr \vu \cdot \Grad \varphi \right] } \dt = - \intO{ \vr_0 \varphi (0, \cdot) }
\]
holds for any test function $\varphi \in C^1_c([0,T) \times \Omega)$;
\item the integral identity
\[
\begin{split}
\int_0^T & \intO{ \left[ \vr \vu \cdot \partial_t \bfphi + \vr \vu \otimes \vu : \Grad \bfphi +
\left( p_{{\rm fl}}(\vr) + p_{{\rm el}}(b) + \frac{2}{3} \sigma b \Del b \right) \Div \bfphi \right] } \dt \\
&= \int_0^T \intO{ \mathbb{S} (\Grad \vu) : \Grad \bfphi } \dt - \sigma
\int_0^T \intO{ \left( \Grad b \otimes \Grad b - \frac{1}{2} |\Grad b|^2 \mathbb{I} \right) : \Grad \bfphi } \dt\\
&- \intO{ \vc{m}_0 \cdot \bfphi (0, \cdot) }
\end{split}
\]
holds for any test function $\bfphi \in C^1_c([0,T) \times \Omega; R^N)$;

\item
\[
\partial_t b + \vu \cdot \Grad b + \frac{1}{\nu}\left( e'(b)  - \sigma \Del b \right)  = \frac{2}{3} b \Div \vu
\]
holds a.a. in $(0,T) \times \Omega$, $b(0, \cdot) = b_0$ a.a. in $\Omega$;

\item the energy inequality
\[
\begin{split}
&\intO{ \left[ \frac{1}{2} \vr |\vu|^2 + \Psi(\vr) + \frac{\sigma}{2} |\Grad b|^2 +  e(b) \right](\tau, \cdot) }\\
&+ \int_0^\tau \intO{ \mathbb{S}(\Grad \vu) : \Grad \vu } + \frac{1}{\nu} \intO{ \left( e'(b) - \sigma \Del b  \right)^2 } \dt \leq E_0
\end{split}
\]
is satisfied for a.a. $\tau \in [0,T]$. Here, $E_0:= \int_{\Omega} \frac{|\vm_0|^2}{2\vr_0} + \Psi(\vr_0) + \frac{\sigma}{2} |\Grad b_0|^2 + e(b_0) \dx$.

\end{itemize}

\end{Definition}

To simplify presentation, we consider hereafter the iconic example:
\begin{equation} \label{A5}
p_{\rm fl}(\vr) = a \vr^\gamma,\
\Psi(\vr) = \frac{a_0}{\gamma - 1} \vr^\gamma, \ e(b) = a_1 b^\alpha - a_2 \log b, \ a_0, a_1, a_2 > 0, \ \alpha \geq 1.
\end{equation}

We are ready to state our main result concerning global--in--time existence of dissipative weak solutions.

\begin{Theorem} \label{Tmain}

Let hypothesis \eqref{A5} be satisfied with
\begin{equation} \label{A9}
\gamma > 3 \textrm{ and } \alpha\ge 1.
\end{equation}
Let the initial data belong to the class
\[
\vr_0 \in L^\gamma(\Omega), \ \vr_0 \geq 0,\ \intO{ \frac{|\vc{m}_0|^2 }{\vr_0} } < \infty,\ b_0 \in W^{1,2}(\Omega) \cap L^\alpha (\Omega),\
\log b_0 \in L^1(\Omega).
\]
Then problem (\ref{M1}--\ref{IC}) admits a dissipative weak solution $[\vr, \vu, b]$ in the sense specified in Definition \ref{D1}.

\end{Theorem}

The remaining part of the paper is devoted to the proof of Theorem \ref{Tmain}. We follow the nowadays standard procedure advocated in the classical monograph by J.-L. Lions
\cite{LIO}, specifically, we establish:

\begin{itemize}
\item {{\it a priori} bounds;}
\item {sequential stability (compactness) of a family o exact solutions;}
\item {suitable approximation scheme.}

\end{itemize}

{
Finally, we note that the assumption $\gamma > 3$ has been made to simplify the proof. It could be possibly relaxed to the nowadays optimal restriction $\gamma > \frac{N}{2}$
by adapting the method developed in \cite{EF70}.
}

\section{Energy estimates and {\it a priori} bounds}
\label{AP}

As mentioned in the preceding section (see also the final part of Appendix),
problem (\ref{M1}--\ref{M6}) admits the total energy of the form
\[
E(\vr, b, \Grad b,  \vu) = \frac{1}{2} \vr |\vu|^2 + \frac{\sigma}{2} |\Grad b |^2 + \Psi(\vr) +  e(b) ,\
\qquad \Psi ' (\vr) \vr - \Psi(\vr) = p_{\rm fl}(\vr).
\]

Indeed, seeing that
\[
\Div \left( \Grad b \otimes \Grad b - \frac{1}{2} |\Grad b |^2 \mathbb{I}\right) =  \Del b \ \Grad b,
\]
we may take the scalar product of (\ref{M2}) with $\vu$, and, integrating the resulting expression over $\Omega$, we obtain
\begin{equation} \label{A1}
\begin{split}
\frac{{\rm d}}{{\rm d}t} &\intO{ \left[ \frac{1}{2} \vr |\vu|^2 + \Psi(\vr) \right] }
+ \intO{ \mathbb{S}(\Grad \vu) :\Grad \vu } \\ &= \intO{ \left( p_{\rm el}(b) + \frac{2}{3} \sigma
b \Del b \right) \Div \vu } - \intO{\sigma \Del b \Grad b  \cdot \vu}
\end{split}
\end{equation}
Next,
multiplying (\ref{M3}) on
\[
 \left( e'(b) - \sigma \Del b \right)
\]
we get
\begin{equation} \label{A2}
\begin{split}
\frac{{\rm d}}{{\rm d}t} &\intO{ \left[ \frac{\sigma}{2} |\Grad b|^2 +  e(b)  \right] }
+  \frac{1}{\nu} \intO{ \left(  e'(b) - \sigma \Del b \right)^2 } \\ &=
- \intO{ \frac{2}{3} \sigma b \Del b \ \Div \vu } +
\intO{  \left( e(b)+\frac{2}{3} b e'(b) \right) \Div \vu } + \sigma
\intO{\vu \cdot \Grad b \ \Del b}.
\end{split}
\end{equation}

Summing up (\ref{A1}), (\ref{A2}) and using relation \eqref{M5}  we may infer that
\begin{equation} \label{A3}
\begin{split}
\frac{{\rm d}}{{\rm d}t} &\intO{ \left[ \frac{1}{2} \vr |\vu|^2 + \Psi(\vr) + \frac{\sigma}{2} |\Grad b|^2 +  e(b) \right] } \\
&+ \intO{ \mathbb{S}(\Grad \vu) : \Grad \vu } + \frac{1}{\nu} \intO{ \left( e'(b) - \sigma \Del b  \right)^2 } = 0.
\end{split}
\end{equation}

In addition to \eqref{A3}, the total mass of the fluid is conserved:
\begin{equation} \label{A4}
\intO{ \vr(t, \cdot) } = \intO{ \vr(0, \cdot) } \ \mbox{for all}\ t \in [0,T].
\end{equation}

To deduce specific estimates from \eqref{A3}, \eqref{A4}, {we use hypothesis (\ref{A5}).}
Consequently, in a accordance with the energy balance (\ref{A3}), we deduce  {\it a priori} bounds in the following function spaces:
\begin{equation} \label{A6}
\begin{split}
\vr &\in L^\infty(0,T; L^\gamma (\Omega)), \\
\vr \vu &\in L^\infty(0,T; L^{\frac{2 \gamma}{\gamma + 1}}(\Omega; R^N)),\\
b &\in L^\infty(0,T; L^\alpha (\Omega) \cap W^{1,2}(\Omega)), \ \log b \in L^\infty (0,T; L^1(\Omega)),\\
\Grad \vu &\in L^2(0,T; L^2(\Omega; R^{N \times N})),\\
- \sigma \Del b + e'(b) &= - \sigma \Del b + a_1 \alpha b^{\alpha - 1} - a_2 \frac{1}{b}
\in L^2((0,T) \times \Omega).
\end{split}
\end{equation}

In addition, we may use (\ref{A4}) and  (\ref{A6}) in combination with a  generalized Poincar\'{e} inequality (see e.g.
\cite[Chapter 11, Section 11.9,11.10]{FEINOV2}) to obtain
\begin{equation} \label{A7}
\vu \in L^2(0,T; W^{1,2}(\Omega; R^N)).
\end{equation}
This together with the first line in \eqref{A6} leads to (by H\"{o}lder's inequality)
\begin{equation}
\vr \vu \in L^2(0,T; L^{\frac{6 \gamma}{\gamma + 6}}(\Omega; R^N)),\label{A7u}
\end{equation}
which, by interpolating with the information given in second line of \eqref{A6}, results at
\begin{equation}
\vr \vu \in L^{\frac{2(5\gamma -3)}{3(\gamma + 1)}}(0,T; L^{\frac{2(5\gamma -3)}{3(\gamma + 1)}}(\Omega; R^N)).\label{A7n}
\end{equation}

Next, seeing that
\[
- \sigma \Del b + e'(b) = h \in L^2((0,T) \times \Omega)
\]
we obtain
\[
\sigma \intO{ |\Del b|^2 } + \intO{ e''(b) |\Grad b|^2 } = - \intO{ h \Del b };
\]
whence, as $e''(b) = a_1 \alpha(\alpha-1) b^{\alpha-2} + \tfrac{a_2}{b^2} \ge \min\{a_1 \alpha(\alpha-1) ,a_2\} \ge 0$,
\begin{equation} \label{A8}
\begin{split}
b &\in L^2(0,T; W^{2,2}(\Omega)), \\
b^{2(\alpha-1)}, \ \frac{1}{b^2} &\in L^1((0,T)\times \Omega),  \\
\Grad b^{\frac{\alpha}{2}}, \ \Grad \log b &\in L^2(0,T; L^2(\Omega; R^N)).
\end{split}
\end{equation}
It follows from \eqref{A6} and \eqref{A8} that, for any  $\alpha\ge 1$,
\begin{equation*}
\begin{split}
b^{\frac{\alpha}{2}} &\in L^{\infty}(0,T; L^{2}(\Omega)) \cap L^2(0,T; W^{1,2}(\Omega)), \\
\log b &\in L^{\infty}(0,T; L^{1}(\Omega)) \cap L^2(0,T; W^{1,2}(\Omega)).
\end{split}
\end{equation*}
Then, by interpolation, we conclude that
\begin{equation}
\begin{split}
b^{\frac{\alpha}{2}} &\in L^{\frac{10}{3}}(0,T; L^{\frac{10}{3}}(\Omega)), \\
\log b &\in L^{\frac{8}{3}}(0,T; L^{\frac{8}{3}}(\Omega)).\end{split} \label{A9a}
\end{equation}
It also follows from \eqref{A6} and \eqref{A9a} that  $b \in L^{\infty}(0,T; W^{1,2}(\Omega)) \cap L^2(0,T; W^{2,2}(\Omega))$. Then, by interpolation and by the embedding theorem, we conclude that
\begin{equation}
\Grad b \in L^{\frac{10}{3}}(0,T; L^{\frac{10}{3}}(\Omega; R^N)) \quad \textrm{ and } \quad b \in L^{20}((0,T)\times \Omega). \label{A9n}
\end{equation}


{The energy estimates established above are strong enough to control all nonlinear terms in the field equations in the reflexive spaces $L^q$, $q > 1$ with the exception of the fluid pressure $p_{{\rm fl}}$. Similarly to the theory of Newtonian compressible fluids, refined estimates are needed to control the pressure team. This is done in the following section.}

\section{The effective viscous flux identity and pressure estimates}
\label{PP}

The effective viscous flux identity is the key tool for proving (i) {\it a priori} bounds yielding the pressure $p_{{\rm fl}}$ in the reflexive space
$L^q$, $q > 1$, (ii) pointwise sequential stability of the density - compactness of $\vr$ in the strong topology of $L^1$. {Here and hereafter, we exploit
essentially the restriction (\ref{A9}), namely}
\[
\gamma > 3.
\]

Let us take the scalar product of the momentum equation \eqref{M2} with the quantity
\[
\Grad \Del^{-1} \left[ \vr - \frac{1}{|\Omega|} \intO{ \vr } \right],
\]
where $\Del^{-1}$ denotes the inverse of the Laplace operator on the space of periodic functions with zero mean. Using the identity
\[
\begin{split}
&\intO{ \Div \mathbb{S}(\Grad \vu) : \Grad \Del^{-1} \left[ \vr - \frac{1}{|\Omega|} \intO{ \vr } \right]} \\ &=
\intO{ \left[ \mu \Del \vu + \left( \frac{1}{3} \mu + \lambda \right) \Grad \Div \vu \right]: \Grad \Del^{-1} \left[ \vr - \frac{1}{|\Omega|} \intO{ \vr } \right]} \\&= - \intO{ \left( \frac{4}{3} \mu + \lambda \right) \vr \Div \vu } ,
\end{split}
\]
we obtain, replacing $\partial_{t} \vr$ by $-\Div (\vr \vu)$,
\begin{equation} \label{A10}
\begin{split}
\int_0^T &\intO{ a_0 \vr^\gamma  \left[ \vr - \frac{1}{|\Omega|} \intO{ \vr } \right]
} \dt \\
&=\left[ \intO{ \vr \vu \cdot \Grad \Del^{-1} \left[ \vr - \frac{1}{|\Omega|} \intO{ \vr } \right] } \right]_{t = 0}^{t=T}
+ \int_0^T \intO{ \vr \vu \cdot \Grad \Del^{-1} [ \Div(\vr \vu)] } \dt \\ &-
\int_0^T \intO{ \vr \vu \otimes \vu : \Grad \Del^{-1} \Grad [ \vr ] } \dt - \int_0^T \intO{ p_{\rm el}(b) \left[ \vr - \frac{1}{|\Omega|} \intO{ \vr } \right]  } \dt  \\
&-\int_0^T \intO{ \frac{2}{3} \sigma b \Del b \left[ \vr - \frac{1}{|\Omega|} \intO{ \vr } \right]  } \dt +
\int_0^T \intO{ \left( \frac{4}{3} \mu + \lambda \right) \vr \Div \vu   } \dt\\
&- \int_0^T \intO{ \sigma \left( \Grad b \otimes \Grad b - \frac{1}{2} |\Grad b |^2 \mathbb{I} \right) : \Grad \Grad \Del^{-1}\left[ \vr  \right]  } = \sum_{j=1}^7 I_j.
\end{split}
\end{equation}
It turns out that all integrals on the right--hand side are controlled in terms of the bounds \eqref{A6}, \eqref{A7}, \eqref{A8} as long as $\gamma > 3$. The restriction $\gamma>3$ comes from estimating $I_2$: noticing that $I_2$ behaves as $\int_0^T \intO{|\vr\vu|^2}$ and referring to \eqref{A7n} one requires that $\frac{2(5\gamma-3)}{3(\gamma+1)} >2$, which gives $\gamma>3$.
As this is well--known for the compressible Navier--Stokes system (see for example \cite{FeiKaPok}), we restrict ourselves to estimating the integrals $I_4$, $I_5$, and $I_7$ containing $b$.

It follows from the hypotheses (\ref{M5}) and \eqref{A5}, (\ref{A6}) and (\ref{A8}) that
\[
p_{\rm el}(b) \approx b^\alpha + \log b  + 1 \in L^{\frac53}((0,T) \times \Omega).
\]
Consequently, the integral $I_4$ remains bounded by the data. Similarly, by virtue of the embedding $W^{1,2} \hookrightarrow L^6$,
\[
b \Del b \in L^2(0,T; L^{\frac{3}{2}}(\Omega));
\]
whence, as $\gamma > 3$, the integral $I_5$ is controlled.
Finally, in accordance with \eqref{A8}, we have
\[
\Grad b \in L^2(0,T; L^6(\Omega; R^N)),
\]
which again yield boundedness of $I_7$ as long as $\gamma > \frac{3}{2}$. Note that
\[
\| \Grad \Grad \Del^{-1}[\vr] \|_{L^q(\Omega)} \approx \left\| \vr - \frac{1}{|\Omega|} \intO{ \vr } \right\|_{L^q(\Omega)}
\ \mbox{for any}\ 1 < q < \infty.
\]
Summing up the previous observations we may deduce from \eqref{A10} that
\begin{equation} \label{A11}
\vr \in L^{\gamma + 1}((0,T) \times \Omega) \ \mbox{if}\ \gamma > 3.
\end{equation}

\section{Sequential stability (compactness)}
\label{S}

At this stage we give a formal argument for showing the weak sequential stability (compactness) of a sequence of
solution $[\vre, \vue, \be]$ of problem (\ref{M1}--\ref{M5}). Similar treatment will be applied to the family of approximate solutions
{introduced in Section \ref{E}}. We merely focus on the steps that are different from the analysis of the standard barotropic Navier--Stokes system referring to \cite{LI4}, \cite{EF70} and \cite{FeiKaPok} for details.

Suppose that $\{ \vre, \vue, \be \}_{\ep > 0}$ is a sequence of solutions to problem (\ref{M1}--\ref{M5}) emanating from the initial data
\begin{equation*}
\begin{split}
\vre(0, \cdot) &= \vr_{0,\ep}, \ \vre \vue (0, \cdot) = \vr_{0, \ep} \vc{v}_{0, \ep}, \ \be(0, \cdot) = b_{0, \ep}, \\
\vr_{0, \ep}, b_{0, \ep} &> 0, \ \vr_{0, \ep} \to \vr_0 \ \mbox{in}\ L^\gamma(\Omega),\ \sup_{\ep>0} \left\{\| \vu_{0,\ep} \|_{L^2(\Omega)} + \| b_{0,\ep} \|_{W^{1,2}(\Omega)} \right\} < \infty.
\end{split}
\end{equation*}
As our arguments will be formal for the time being, we suppose we deal with strong solutions, however, as we shall see below, the same treatment can be applied
to families of weak solutions. {To simplify presentation, we also consider the initial data more regular than in Theorem \ref{Tmain}, specifically,
$\vr_{0, \ep}$ and $b_{0, \ep}$ are supposed to be bounded below away from zero uniformly for $\ep \to 0$.}

In accordance with the {\it a priori} bounds established in Sections \ref{AP}, \ref{PP}, the estimates (\ref{A6}--\ref{A8}),
\eqref{A11} hold uniformly for $\ep \to 0$. Passing to suitable subsequences, we may therefore assume that
\[
\begin{split}
\vre &\to \vr \ \mbox{in}\ C_{\rm weak}([0,T]; L^\gamma(\Omega)) \ \mbox{and weakly in}\ L^{\gamma + 1}((0,T) \times \Omega),\\
\vue &\to \vu \ \mbox{weakly in}\ L^2(0,T; W^{1,2}_0 (\Omega; R^N)), \\
\be &\to b \ \mbox{weakly in}\ L^2(0,T; W^{2,2}(\Omega)) \ \mbox{and in}\ C([0,T]; W^{1,2}(\Omega)).
\end{split}
\]
In addition, as $\be$ satisfies the standard parabolic equation
\[
\partial_t \be - \frac{\sigma}{\nu} \Del \be = - \vu \cdot \Grad \be - \frac{e'(\be)}{\nu} + \sigma \frac{2}{3} \be \Div \vue,
\]
where the right--hand side belongs to $L^q((0,T) \times \Omega)$ for a certain $q > 1$ (cf. \eqref{A6}, \eqref{A8}), we get
\begin{equation} \label{A12}
\begin{split}
\partial_t \be &\to \partial_t b , \ \Del \be \to \Del b \ \mbox{weakly in}\ L^q((0,T) \times \Omega),\\
\mbox{in particular}\ \be &\to b,\ \Grad \be \to \Grad b \ \mbox{(strongly) in}\ L^q((0,T) \times \Omega; R^N).
\end{split}
\end{equation}

Now, following the arguments of \cite[Chapter 6]{FeiKaPok}, we can let $\ep \to 0$ in (\ref{M1}--\ref{M3}) to obtain
\begin{equation} \label{A13}
\partial_t \vr + \Div (\vr \vu) = 0,
\end{equation}
\begin{equation} \label{A14}
\begin{split}
\partial_t (\vr \vu) + \Div (\vr \vu \otimes \vu) &+ \Grad \left( a \Ov{\vr^\gamma} + p_{{\rm el}}(b)  + \frac{2}{3}
\sigma b \Del b \right) \\
&= \Div \mathbb{S} (\Grad \vu) - \sigma \Div \left( \Grad b \otimes \Grad b - \frac{1}{2} |\Grad b |^2\right),
\end{split}
\end{equation}
\begin{equation} \label{A15}
\partial_t b + \vu \cdot \Grad b + \frac{1}{\nu}\left( e'(b)  - \sigma \Del b\right) = \frac{2}{3} b \Div \vu,
\end{equation}
where the equations (\ref{A13}), (\ref{A14}) are satisfied in the sense of distributions, while \eqref{A15} holds a.a. in $(0,T) \times \Omega$.
The symbol $\Ov{\vr^\gamma}$ denotes a weak limit of $\vre^\gamma$. Consequently, it remains to show
\[
\Ov{\vr^\gamma} = \vr^\gamma \ \mbox{or, equivalently, strong convergence}\ \vre \to \vr \ \mbox{in}\  L^1((0,T) \times \Omega).
\]

Evoking the method developed by Lions \cite{LI4}, we first observe that $\vr$, $\vu$ satisfy the renormalized equation of continuity,
\begin{equation} \label{A16}
\partial_t \beta(\vr) + \Div (\beta(\vr) \vu) + \left( \beta'(\vr) \vr - \beta(\vr) \right) \Div \vu = 0
\end{equation}
for any $\beta$ with a sufficiently moderate growth for large $\vr$. Indeed, as $\gamma > 3$, equation (\ref{A16}) follows from (\ref{A13}) by means of the regularization procedure proposed by DiPerna and Lions \cite{DL}. In particular, using compactness of the initial densities, we deduce from (\ref{A16}) and the corresponding equation for $\beta_{\varepsilon}(\vr)$ that
\begin{equation} \label{A17}
\intO{ \left[ \Ov{ \vr \log \vr } - \vr \log \vr  \right](\tau, \cdot) } + \int_0^\tau \intO{ \left( \Ov{ \vr \Div \vu} -
\vr \Div \vu \right) } \dt  = 0 \ \mbox{for any}\ \tau \in [0,T].
\end{equation}
As the function $\vr \mapsto \vr \log \vr $ is strictly convex, the desired strong convergence of $\{ \vre \}_{\ep > 0}$ follows as soon as we show
\begin{equation} \label{A18}
\int_0^\tau \intO{ \left( \Ov{ \vr \Div \vu} -
\vr \Div \vu \right) } \dt \geq 0 \ \mbox{for any}\ \tau \geq 0.
\end{equation}

To see (\ref{A18}) we use the celebrated \emph{effective viscous flux identity} of Lions \cite{LI4}. We start with the relation (\ref{A10}) evaluated at the level of $\ep$-approximation:
\begin{equation} \label{A19}
\begin{split}
\int_0^\tau &\intO{ a_0 \vre^\gamma  \left[ \vre - \frac{1}{|\Omega|} \intO{ \vre } \right]
} \dt \\
&=\left[ \intO{ \vre \vue \cdot \Grad \Del^{-1} \left[ \vre - \frac{1}{|\Omega|} \intO{ \vre } \right] } \right]_{t = 0}^{t=\tau}
+ \int_0^\tau \intO{ \vre \vue \cdot \Grad \Del^{-1} [ \Div(\vre \vue) ] } \dt \\ &-
\int_0^\tau \intO{ \vre \vue \otimes \vue : \Grad \Del^{-1} \Grad [ \vre ] } \dt - \int_0^\tau \intO{ p_{\rm el}(\be) \left[ \vre - \frac{1}{|\Omega|} \intO{ \vre } \right]  } \dt  \\
&-\int_0^\tau \intO{ \frac{2}{3} \sigma \be \Del b \left[ \vre - \frac{1}{|\Omega|} \intO{ \vre } \right]  } \dt +
\int_0^\tau \intO{ \left( \frac{4}{3} \mu + \lambda \right) \vre \Div \vue   } \dt\\
&- \int_0^\tau \intO{ \sigma \left( \Grad \be \otimes \Grad \be - \frac{1}{2} |\Grad \be |^2 \mathbb{I}\right) : \Grad \Grad \Del^{-1}\left[ \vre  \right]  } \dt .
\end{split}
\end{equation}

Similarly, we may deduce from the limit system (\ref{A13}--\ref{A15}) the identity
\begin{equation} \label{A20}
\begin{split}
\int_0^\tau &\intO{ a_0 \Ov{\vr^\gamma}  \left[ \vr - \frac{1}{|\Omega|} \intO{ \vr } \right]
} \dt \\
&=\left[ \intO{ \vr \vu \cdot \Grad \Del^{-1} \left[ \vr - \frac{1}{|\Omega|} \intO{ \vr } \right] } \right]_{t = 0}^{t=\tau}
+ \int_0^\tau \intO{ \vr \vu \cdot \Grad \Del^{-1} [ \Div(\vr \vu)] } \dt \\ &-
\int_0^\tau \intO{ \vr \vu \otimes \vu : \Grad \Del^{-1} \Grad [ \vr ] } \dt - \int_0^\tau \intO{ p_{\rm el}(b) \left[ \vr - \frac{1}{|\Omega|} \intO{ \vr } \right]  } \dt  \\
&-\int_0^\tau \intO{ \frac{2}{3} \sigma b \Del b \left[ \vr - \frac{1}{|\Omega|} \intO{ \vr } \right]  } \dt +
\int_0^\tau \intO{ \left( \frac{4}{3} \mu + \lambda \right) \vr \Div \vu   } \dt\\
&- \int_0^\tau \intO{ \sigma \left( \Grad b \otimes \Grad b - \frac{1}{2} |\Grad b |^2 \mathbb{I} \right) : \Grad \Grad \Del^{-1}\left[ \vr  \right]  } \dt .
\end{split}
\end{equation}

Now, using the compactness arguments known for the compressible Navier--Stokes system (cf. e.g. \cite[Chapter 6]{FeiKaPok}), we may let $\ep \to 0$
in (\ref{A19}) and compare
the resulting expression with (\ref{A20}) obtaining:
\begin{equation} \label{A21}
\begin{split}
0 \leq
\lim_{\ep \to 0} &\int_0^\tau \intO{ a \vre^{\gamma + 1} } \dt -
\int_0^\tau \intO{ a \Ov{\vr^\gamma} \vr
} \dt \\
&=\int_0^\tau \intO{ \frac{2}{3} \sigma b \Big( \vr \Del b  - \Ov{ \vr \Del b} \Big)  } \dt \\ &+
\int_0^\tau \intO{ \left( \frac{4}{3} \mu + \lambda \right)\Big( \Ov{\vr \Div \vu} - \vr \Div \vu \Big)   } \dt.
\end{split}
\end{equation}
In comparison with the analysis of the standard Navier--Stokes system, we have to handle the first integral on the right--hand side of \eqref{A21}.
To continue, we recall that (\ref{A21}) can be localized (cf. \cite[Chapter 6]{FeiKaPok}), namely
\[
\begin{split}
0 \leq
\lim_{\ep \to 0} &\int_0^T \intO{ a \vre^{\gamma + 1} \varphi } \dt -
\int_0^T \intO{ a \Ov{\vr^\gamma} \vr \varphi
} \dt \\
&=\int_0^T \intO{ \frac{2}{3} \sigma b \Big( \vr \Del b  - \Ov{ \vr \Del b} \Big) \varphi  } \dt \\ &+
\int_0^T \intO{ \left( \frac{4}{3} \mu + \lambda \right)\Big( \Ov{\vr \Div \vu} - \vr \Div \vu \Big) \varphi   } \dt.
\end{split}
\]
for any $\varphi \in \DC((0,T) \times \Omega)$, $\varphi \geq 0$. In other words, almost everywhere in $(0,T)\times \Omega$, we have
\begin{equation} \label{A21a}
\frac{2}{3} \sigma b \Big( \vr \Del b  - \Ov{ \vr \Del b} \Big) + \left( \frac{4}{3} \mu + \lambda \right)\Big( \Ov{\vr \Div \vu} - \vr \Div \vu \Big)
\geq 0.
\end{equation}

Now,
at the level of $\ep$-approximation, we have
\[
\partial_t \be + \vue \cdot \Grad \be + \frac{1}{\nu} \left(e'(\be)  - \sigma \Del \be\right) = \frac{2}{3} \be \Div \vue,
\]
or
\[
\partial_t (\vre \be) + \Div (\vre \be \vue) + \frac{1}{\nu} \left(\vre e'(\be)  - \sigma \vre \Del \be \right) = \frac{2}{3} \be \vre \Div \vue.
\]
Letting $\ep \to 0$ we get
\begin{equation} \label{A22}
\partial_t (\vr b) + \Div (\vr b \vu) + \frac{1}{\nu}\left( \vr e'(b)  - \sigma \Ov{ \vr \Del b } \right) = \frac{2}{3} b \Ov{ \vr \Div \vu}
\end{equation}
in the sense of distributions.

On the other hand, multiplying the limit equation (\ref{A15}) on $\vr$ and using (\ref{A13}), we obtain
\begin{equation} \label{A23}
\partial_t (\vr b) + \Div (\vr b \vu) + \frac{1}{\nu} \left( \vr e'(b)  - \sigma \vr \Del b \right)  = \frac{2}{3} b \vr \Div \vu.
\end{equation}
Thus comparing (\ref{A22}), (\ref{A23}), we deduce
\[
\vr \Del b - \Ov{\vr \Del b} = \frac{2}{3 \nu} b \Big( \Ov{\vr \Div \vu} - \vr \Div \vu \Big),
\]
which, together with (\ref{A21a}) gives rise to
\[
\left( \frac{4 \sigma b^2}{9\nu} + \frac{4}{3} \mu + \lambda \right) \Big( \Ov{\vr \Div \vu} - \vr \Div \vu \Big)  \geq 0.
\]
Thus we have shown (\ref{A18}), and, consequently, the desired strong convergence
\[
\vre \to \vr \ \mbox{in}\ L^1((0,T) \times \Omega).
\]

\section{Approximation scheme}
\label{E}

To construct the weak solutions to problem (\ref{M1}--\ref{M3}) we adapt the multilevel approximation scheme from \cite[Chapter 7]{FeiKaPok}:
\begin{equation} \label{E1}
\partial_t \vr + \Div (\vr \vu) = \ep \Del \vr,
\end{equation}
\begin{equation} \label{E2}
\begin{split}
\partial_t \Pi_m (\vr \vu) &+ \Div \Pi_m (\vr \vu \otimes \vu) + \Grad \Pi_m \left( p_{\rm fl}(\vr) + p_{\rm el}(b) + \frac{2}{3} \sigma b \Del b \right) \\ &=
\ep \Del \Pi_m (\vr \vu) + \Div \Pi_m \mathbb{S} (\Grad \vu) -
\sigma \Div \Pi_m \left( \Grad b \otimes \Grad b - \frac{1}{2} |\Grad b|^2 \mathbb{I} \right),
\end{split}
\end{equation}
\begin{equation} \label{E3}
\partial_t b + \vu \cdot \Grad b + \frac{1}{\nu} \left( e'(b)  - \sigma \Del b \right) = \frac{2}{3} b \Div \vu.
\end{equation}
Here, $\Pi_m$ denotes the orthogonal projection on a finite dimensional space $X_m$ spanned by a finite number of trigonometric polynomials.

\subsection{Solvability of the Galerkin approximation}

For fixed $m > 0$, $\ep > 0$, the approximate problem (\ref{E1}--\ref{E3}) can be solved by a fixed point argument, similarly to \cite[Chapter 7]{FeiKaPok}.
More specifically, given $b \in C([0,T]; W^{1,2}(\Omega))$, the system (\ref{E1}), (\ref{E2}) admits a unique solution $(\vr, \vu) =
(\vr, \vu)[b]$. Plugging $\vu = \vu[b]$ in (\ref{E3}), we can find a unique solution $\mathcal{T}[b]$ of (\ref{E3}). Given compactness properties
of solutions to the parabolic equation (\ref{E3}), the solution of (\ref{E1}--\ref{E3}) can be obtained as a fixed point of the mapping $\mathcal{T}$.
We recall that all the {\it a priori} bounds derived formally in Section \ref{AP} remain valid at this level of approximation as they are based
the energy estimates that are compatible with the Galerkin approximation.

Exactly as in \cite[Chapter 7]{FeiKaPok}, we can pass to the limit $m \to \infty$ obtaining the second level approximate solutions satisfying
\begin{equation} \label{E4}
\partial_t \vr + \Div (\vr \vu) = \ep \Del \vr,
\end{equation}
\begin{equation} \label{E5}
\begin{split}
&\left[ \intO{ \vr \vu \cdot \bfphi } \right]_{t = 0}^{t = \tau} =
\int_0^\tau \intO{ \left[ \vr \vu \cdot \partial_t \bfphi + \vr \vu \otimes \vu : \Grad \bfphi  + \left( p_{\rm fl}(\vr) + p_{\rm el}(b) + \frac{2}{3} \sigma b \Del b \right) \Div \bfphi \right] } \dt  \\ &= \int_0^\tau \intO{ \left[
\ep (\vr \vu) \cdot \Del \bfphi - \mathbb{S} (\Grad \vu) : \Grad \bfphi +
\sigma \left( \Grad b \otimes \Grad b - \frac{1}{2} |\Grad b|^2 \mathbb{I} \right) : \Grad \bfphi \right]} \dt,
\end{split}
\end{equation}
\begin{equation} \label{E6}
\partial_t b + \vu \cdot \Grad b + \frac{1}{\nu}\left(e'(b)  - \sigma \Del b\right) = \frac{2}{3} b \Div \vu.
\end{equation}

\subsection{Vanishing viscosity limit}

Our ultimate goal is to let $\ep \to 0$ in the approximate problem (\ref{E4}--\ref{E5}). The reader will notice that this step differs from the formal
procedure performed in Section \ref{S} only by the fact that $\be$ satisfies the perturbed equation
\[
\partial_t (\vre \be) + \Div (\vre \vue \be) + \vre e'(\be)  - \sigma \vre \Del \be = \frac{2}{3} \vre \be \Div \vue -
\ep \be \Del \vre,
\]
where the extra $\ep-$term vanishes in the asymptotic limit.

{Having performed this last step, we have completed the proof of Theorem \ref{Tmain}.}

\section{Appendix}
\label{Z}

In this section we derive the model analyzed in this paper. We use a methodology developed within the context of incompresssible viscoelastic fluids in \cite{RASR}. Here, we follow more recent studies concerning compressible viscoelastic rate-type fluids without or with stress diffusion presented in \cite{mp2016} and in \cite{mpss2018, bmps2018}, where the reader can found more detailed expositions.

The methodology is based on the tenet that as the body dissipates the energy (or more generally produce the entropy) there is an underlying evolving natural configuration associated with the current configuration and the response between these configurations is purely elastic. The natural configuration thus splits the total deformation process into the part that is elastic and the other part that describes the (irreversible) changes in the natural configuration. The governing constitutive equations for the Cauchy stress are determined from the knowledge of the constitutive equations for two scalar quantities: the Helmholtz free energy and the rate of entropy production. The constitutive equation for the Helmholtz free energy specifies (requires to know) how the material stores the energy, while the constitutive equation for the rate of entropy production describes how the material dissipates the energy.

Let $\vv$ be the velocity defined at the current configuration $\kappa_t$ and $\mathbb{F}$ be the deformation tensor acting between the reference configuration $\kappa_R$ and the configuration $\kappa_t$. These quantities are linked through the equation for the material time derivative of $\mathbb{F}$ that takes the form
\begin{equation}
\tecka{\mathbb{F}} = \Grad \vv \, \mathbb{F} \qquad \textrm{ or } \qquad \Grad \vv = \tecka{\mathbb{F}} \mathbb{F}^{-1}. \label{aA1}
\end{equation}
We also define
\begin{equation}
  \mathbb{B}:= \mathbb{F} \mathbb{F}^{T}, \qquad \mathbb{C}:= \mathbb{F}^{T} \mathbb{F}, \qquad \mathbb{D}:= \frac{\Grad \vv + (\Grad \vv)^{T} }{2}, \label{aA1a}
\end{equation}
where $\mathbb{A}^T$ stands for the transpose matrix of $\mathbb{A}$.

The natural configuration $\kappa_{p(t)}$, evolving with the current configuration $\kappa_t$, splits the total deformation $\mathbb{F}$ into its elastic part $\mathbb{F}_{\kappa_{p(t)}}$ and the dissipative part $\mathbb{G}$ so that
\begin{equation}
\mathbb{F} = \mathbb{F}_{\kappa_{p(t)}} \mathbb{G}. \label{aA2}
\end{equation}
In virtue of \eqref{aA1} and \eqref{aA1a} we also set
\begin{equation}
\mathbb{L}_{\kappa_{p(t)}} := \tecka{\mathbb{G}} \mathbb{G}^{-1}, \qquad \mathbb{D}_{\kappa_{p(t)}}:= \frac{\mathbb{L}_{\kappa_{p(t)}} + (\mathbb{L}_{\kappa_{p(t)}})^T}{2} \label{aA3}
\end{equation}
and
\begin{equation}
\mathbb{B}_{\kappa_{p(t)}} := \mathbb{F}_{\kappa_{p(t)}} \mathbb{F}_{\kappa_{p(t)}}^T, \qquad \mathbb{C}_{\kappa_{p(t)}} := \mathbb{F}_{\kappa_{p(t)}}^T \mathbb{F}_{\kappa_{p(t)}}.\label{aA4}
\end{equation}
It follows from these relations that
\begin{equation}
\tecka{\mathbb{B}}_{\kappa_{p(t)}} = (\Grad \vv) \mathbb{B}_{\kappa_{p(t)}}  + \mathbb{B}_{\kappa_{p(t)}} (\Grad \vv)^T - 2 \mathbb{F}_{\kappa_{p(t)}} \mathbb{D}_{\kappa_{p(t)}} \mathbb{F}_{\kappa_{p(t)}}^T \label{aA5}
\end{equation}
and
\begin{equation}
\tr \tecka{\mathbb{B}}_{\kappa_{p(t)}} = 2 \mathbb{D}:\mathbb{B}_{\kappa_{p(t)}}  - 2 \mathbb{C}_{\kappa_{p(t)}}: \mathbb{D}_{\kappa_{p(t)}}. \label{aA6}
\end{equation}

We now introduce a key assumption leading to the model analyzed in this study. We will require that $\mathbb{B}_{\kappa_{p(t)}}$ is only spherical. This means that the elastic response is connected merely with the expansion/contraction or in some simple cases with growth/degradation. Hence
\begin{equation}
\mathbb{B}_{\kappa_{p(t)}} = b \, \mathbb{I} \quad \textrm{ where } b := \frac{\tr \mathbb{B}_{\kappa_{p(t)}}}{3}. \label{aA7}
\end{equation}
This implies that $(\mathbb{B}_{\kappa_{p(t)}})_{\delta}$ is zero and $\tr \mathbb{B}_{\kappa_{p(t)}} = \tr \mathbb{C}_{\kappa_{p(t)}}$. Consequently, $(\mathbb{C}_{\kappa_{p(t)}})_{\delta}$ is zero and $\mathbb{C}_{\kappa_{p(t)}} = b \mathbb{I}$.
Inserting \eqref{aA7} into \eqref{aA6} we obtain
\begin{equation}
\tecka{b} = \frac23 b \Div \vv - \frac23 b \tr \mathbb{D}_{\kappa_{p(t)}}. \label{aA8}
\end{equation}

Next, we recall the balance equations for mass, linear and angular momenta. These take the form
\begin{equation}
\tecka{\varrho} = - \varrho \Div \vv,\quad \varrho \tecka{\vv} = \Div \mathbb{T}, \qquad \mathbb{T}=\mathbb{T}^T,  \label{aA9}
\end{equation}
where $\varrho$ is the density and $\mathbb{T}$ is the Cauchy stress.

Restricting ourselves to isothermal processes, the formulation of the first and the second law thermodynamics (the balance of energy and the balance of entropy) reduces to (see \cite{bmps2018} for details)
\begin{equation}
\mathbb{T}:\mathbb{D} - \varrho \tecka{\psi} = \xi - \Div (\bfj_{e} - \theta \bfj_{\eta}) \quad \textrm{ with } \xi \ge 0, \label{aA10}
\end{equation}
where $\psi$ is the Helmholtz free energy and $\xi$ is the rate of dissipation.

The starting point of the constitutive theory is the assumption how the material stores the energy, which characterizes the response between natural and current configuration. In \cite{mpss2018} and \cite{bmps2018}, the authors require, following also the earlier results presented in \cite{RASR}, \cite{mtr2015} and \cite{mp2016}, that
\begin{equation}
\psi = \psi_0(\varrho) + \frac{\mu}{2\varrho} \left( \tr \mathbb{B}_{\kappa_{p(t)}} -3 - \ln \det \mathbb{B}_{\kappa_{p(t)}} \right) + \frac{\sigma}{2\varrho} |\Grad \frac{\tr \mathbb{B}_{\kappa_{p(t)}}}{3} |^2. \label{aA11}
\end{equation}
Considering $\mathbb{B}_{\kappa_{p(t)}}$ of the particular form \eqref{aA7}, the constitutive assumption \eqref{aA11} reduces to
\begin{equation}
\psi = \psi_0(\varrho) + \frac{1}{\varrho} \left( e(b) + \frac{\sigma}{2} |\Grad b|^2 \right) \quad \textrm{ with } e(b) := \frac{3\mu}{2}(b-1- \ln b). \label{aA12}
\end{equation}
This implies, using also \eqref{aA8} and the first equation in \eqref{aA9}, that
\begin{equation}
\begin{split}
 \varrho \tecka{\psi} &= \Div (\sigma \tecka{b} \Grad b) - \left(\varrho^2 \psi_0'(\varrho) - (e(b) + \frac23 b e'(b)) - \frac{\sigma}{2} |\Grad b|^2 + \sigma b \Delta b\right) \Div \vv \\
& - \sigma (\Grad b \otimes \Grad b):\mathbb{D} - \frac23 (b e'(b) - \sigma b \Delta b) \tr \mathbb{D}_{\kappa_{p(t)}}\,.
\end{split} \label{aA13}
\end{equation}
Denoting
\begin{equation}
  p_{\textrm{fl}}(\varrho):= \varrho^2 \psi_0'(\varrho) \quad \textrm{ and } \quad p_{\textrm{el}}(b):= - e(b) - \frac{2}{3} b e'(b), \label{aA000}
\end{equation}
setting
\begin{equation}
  \theta \bfj_{\eta} := \bfj_{e} -  \sigma\tecka{b}\Grad b,\label{aA00}
\end{equation}
and inserting \eqref{aA13} into \eqref{aA10}, we obtain
\begin{equation}
\begin{split}
 \xi &= (\mathbb{T}_{\delta} + \sigma (\Grad b \otimes \Grad b)_{\delta}) : \mathbb{D}_{\delta} \\
&+ (m+ {3} |\Grad b|^2 +  p_{\textrm{fl}}(\varrho)
+ p_{\textrm{el}}(b) - \frac{\sigma}{2} |\Grad b|^2 + \sigma b \Delta b)\Div \vv \\
& + \frac23 (b e'(b) - \sigma b \Delta b) \tr \mathbb{D}_{\kappa_{p(t)}}\,,
\end{split} \label{aA14}
\end{equation}
where we used the decomposition $\mathbb{A}= \mathbb{A}_{\delta} +\frac13 (\tr \mathbb{A}) \mathbb{I}$; $\mathbb{A}_{\delta}$ denotes the deviatoric (traceless) part of $\mathbb{A}$ and $\mathbb{I}$ is the identity tensor. In \eqref{aA14}, $m$ stands for the mean normal stress, i.e. $m:=\dfrac13 \tr \mathbb{T}$.

In order to guarantee apriori that $\xi \ge 0$ we set
\begin{align}
\mathbb{T}_{\delta} + \sigma (\Grad b \otimes \Grad b)_{\delta} &= 2\mu \mathbb{D}_{\delta}, \label{Aa15}\\
m+ {3} |\Grad b|^2 +  p_{\textrm{fl}}(\varrho) + p_{\textrm{el}}(b) - \frac{\sigma}{2} |\Grad b|^2 + \sigma b \Delta b &= \lambda \Div \vv \label{aA16}\\
 (b e'(b) - \sigma b \Delta b) &= \frac23 \nu b^{s} \tr \mathbb{D}_{\kappa_{p(t)}}\,. \label{aA17}
\end{align}
Here $\mu$, $\lambda$ and $\nu$ are positive and $s\ge 0$. Note that it follows from the definition of $e$, see \eqref{aA12}, that $b$ is positive.

It follows from \eqref{Aa15} and \eqref{aA16} that
\begin{align}
\mathbb{T} &= \mathbb{T}_{\delta} + m\mathbb{I} = - (p_{\textrm{fl}}(\varrho)+ p_{\textrm{el}}(b)) \mathbb{I} - \sigma (\Grad b \otimes \Grad b - \frac{1}{2} |\Grad b|^2\mathbb{I}) + 2\mu \mathbb{D}_{\delta} + \lambda \Div \vv \mathbb{I}, \label{aA18}
\end{align}
while \eqref{aA17} and \eqref{aA8} leads to
\begin{align}
\tecka{b} + \frac{1}{\nu} (e'(b) - \sigma \Delta b) b^{2-s} &= \frac23 b \Div \vv \,. \label{aA19}
\end{align}
Upon inserting these equations into \eqref{aA14}, we observe that
\begin{equation}
 \xi = 2\mu |\mathbb{D}_{\delta}|^2 + \lambda (\Div \vv)^2 + \frac{1}{\nu} (e'(b) - \sigma \Delta b)^2 b^{2 - s}.
\label{aA22}
\end{equation}
We also observe that the constitutive equations \eqref{aA18}-\eqref{aA19} coincide with \eqref{M3} if $s=2$.

Taking the scalar product of $\varrho \tecka{\vv} = \Div \mathbb{T}$ and $\vv$, we obtain
\begin{equation*}
\varrho \frac{\partial}{\partial t} (\frac{|\vv|^2}{2}) + \vr \vv \cdot \Grad \frac{|\vv|^2}{2} - \Div (\mathbb{T}\vv) + \mathbb{T}:\mathbb{D} = 0.
\end{equation*}
Using then the balance equation for mass, see \eqref{aA9}$_1$, and the reduced thermodynamical identity \eqref{aA10} together with \eqref{aA00}, we arrive at
\begin{equation}
 \frac{\partial}{\partial t} (\vr \frac{|\vv|^2}{2}) + \Div( \vr \frac{|\vv|^2}{2} \vv - \mathbb{T}\vv - \sigma\tecka{b}\Grad b) + \xi  + \vr\tecka{\psi} = 0, \label{aA23}
\end{equation}
where $\xi$ fulfills \eqref{aA22} and $\psi$ fulfills \eqref{aA12}. Next, using the identity (derived with help of \eqref{aA9}$_1$ again)
\begin{equation*}
\vr \tecka{\psi} = (\vr \psi\tecka{)} - \tecka{\vr} \psi = \frac{\partial}{\partial t} (\vr \psi) + \Div (\vr\psi\vv),
\end{equation*}
we conclude from \eqref{aA23} that
\begin{equation}
 \frac{\partial}{\partial t} (\frac12 \vr |\vv|^2 + \vr \psi) + \Div(\frac{1}{2} \vr |\vv|^2 \vv  + \vr \psi \vv - \mathbb{T}\vv - \sigma\tecka{b}\Grad b) + \xi  = 0. \label{aA24}
\end{equation}
Integrating this identity over $\Omega$, using the Gauss theorem and assuming that the boundary integrals vanish, we get
\begin{equation}
 \frac{d}{dt} \int_{\Omega} ( \frac12 \vr |\vv|^2  + \vr \psi)  + \int_{\Omega} \xi  = 0. \label{aA25}
\end{equation}
Finally, using \eqref{aA12} and \eqref{aA22}, and setting
\begin{equation}
\Psi(\vr):= \vr\psi_0 (\vr), \label{aA20}
\end{equation}
we conclude from \eqref{aA25} the energy identity
\begin{equation}
\begin{split}
 \frac{d}{dt} \int_{\Omega} &(\frac12 \vr |\vv|^2  + \Psi(\vr) + e(b) + \frac{\sigma}{2} |\Grad b|^2)  \\
&+ \int_{\Omega} 2\mu |\mathbb{D}_{\delta}|^2 + \lambda (\Div \vv)^2 + \frac{1}{\nu} (e'(b) - \sigma \Delta b)^2  = 0.
\end{split}\label{aA26}
\end{equation}

Note that it follows from \eqref{aA20} and \eqref{aA000} that
\begin{equation}
\vr \Psi'(\vr) - \Psi(\vr) = p_{\textrm{fl}}(\vr). \label{aA20c}
\end{equation}

\bigskip

\centerline{\bf Acknowledgement}

\medskip

The authors acknowledge support of the project 18-12719S financed
by the Czech Science Foundation.

\bigskip

\def\cprime{$'$} \def\ocirc#1{\ifmmode\setbox0=\hbox{$#1$}\dimen0=\ht0
  \advance\dimen0 by1pt\rlap{\hbox to\wd0{\hss\raise\dimen0
  \hbox{\hskip.2em$\scriptscriptstyle\circ$}\hss}}#1\else {\accent"17 #1}\fi}


\begin{thebibliography}{10}

\bibitem{AKM}
S.~N. Antontsev, A.~V. Kazhikhov, and V.~N. Monakhov.
\newblock {\em Krajevyje zadaci mechaniki neodnorodnych zidkostej}.
\newblock Novosibirsk, 1983.

\bibitem{bmps2018}
M.~Bul\'\i\v{c}ek, J.~M\'alek, V.~Pr\ocirc{u}\v{s}a, and E.~S\"uli.
\newblock P{DE} analysis of a class of thermodynamically compatible
  viscoelastic rate-type fluids with stress-diffusion.
\newblock In {\em Mathematical analysis in fluid mechanics---selected recent
  results}, volume 710 of {\em Contemp. Math.}, pages 25--51. Amer. Math. Soc.,
  Providence, RI, 2018.

\bibitem{DL}
R.J. DiPerna and P.-L. Lions.
\newblock Ordinary differential equations, transport theory and {S}obolev
  spaces.
\newblock {\em Invent. Math.}, {\bf 98}:511--547, 1989.

\bibitem{EF70}
E.~Feireisl.
\newblock {\em Dynamics of viscous compressible fluids}.
\newblock Oxford University Press, Oxford, 2004.

\bibitem{FeJiNo}
E.~Feireisl, Bum~Ja Jin, and A.~Novotn{\' y}.
\newblock Relative entropies, suitable weak solutions, and weak-strong
  uniqueness for the compressible {N}avier-{S}tokes system.
\newblock {\em J. Math. Fluid Mech.}, {\bf 14}:712--730, 2012.

\bibitem{FeiKaPok}
E.~Feireisl, T.~Karper, and M.~Pokorn{\' y}.
\newblock {\em Mathematical theory of compressible viscous fluids: {A}nalysis
  and numerics}.
\newblock Birkh{\" a}user--Verlag, Basel, 2017.

\bibitem{FEINOV2}
E.~Feireisl and A.~Novotn\'y.
\newblock {\em Singular limits in thermodynamics of viscous fluids}.
\newblock Advances in Mathematical Fluid Mechanics. Birkh\"auser/Springer,
  Cham, 2017.

\bibitem{LIO}
J.-L. Lions.
\newblock {\em Quelques m{\' e}thodes de r{\' e}solution des probl{\` e}mes aux
  limites non lin{\' e}aires}.
\newblock Dunod, Gautthier - Villars, Paris, 1969.

\bibitem{LI4}
P.-L. Lions.
\newblock {\em Mathematical topics in fluid dynamics, Vol.2, Compressible
  models}.
\newblock Oxford Science Publication, Oxford, 1998.

\bibitem{mpss2018}
J.~M\'{a}lek, V.~Pr\ocirc{u}\v{s}a, T.~Sk\v{r}ivan, and E.~S\"{u}li.
\newblock Thermodynamics of viscoelastic rate-type fluids with stress
  diffusion.
\newblock {\em Physics of Fluids}, 30(2):023101, 2018.

\bibitem{mp2016}
J.~M{\'a}lek and V.~Pr{\r{u}}{\v{s}}a.
\newblock {\em Derivation of Equations for Continuum Mechanics and
  Thermodynamics of Fluids}, pages 1--70.
\newblock Springer International Publishing, Cham, 2016.

\bibitem{mtr2015}
J.~M\'{a}lek, K.R. Rajagopal, and K.~T\ocirc{u}ma.
\newblock On a variant of the {M}axwell and {O}ldroyd-{B} models within the
  context of a thermodynamic basis.
\newblock {\em International Journal of Non-Linear Mechanics}, {\bf 76}:42--47,
  2015.

\bibitem{RASR}
K.R. Rajagopal and A.R. Srinivasa.
\newblock A thermodynamic frame work for rate type fluid models.
\newblock {\em Journal of Non-Newtonian Fluid Mechanics}, 88(3):207 -- 227,
  2000.

\end{thebibliography}

\end{document}